\documentclass[a4paper,11pt]{article}

\usepackage{amsfonts}
\usepackage{amsmath}
\usepackage{amstext}
\usepackage{amsthm}
\usepackage{xspace}

\input psfig.sty

\newcommand{\N}{\mathbb N}

\newcommand{\E}{\mathbb E}

\newtheorem{teorema} {Theorem}[section]
\newtheorem{prop}[teorema]{Proposition}

\theoremstyle{definition}
\newtheorem{definizione}[teorema]{Definition}

\newtheorem{guess}[teorema]{Remark}

\begin{document}
\title{The Algorithmic Information Content for randomly perturbed
  systems} \author{Claudio Bonanno\footnote{Dipartimento di
    Matematica, Universit\`a di Pisa, via Buonarroti 2/a, 56127 Pisa
    (Italy), email: bonanno@dm.unipi.it}} \date{} \maketitle

\begin{abstract}
In this paper we prove estimates on the behaviour of the
Kolmo\-go\-rov-Sinai entropy relative to a partition for randomly
perturbed dynamical systems. Our estimates use the entropy for the
unperturbed system and are obtained using the notion of Algorithmic
Information Content. The main result is an extension of known results
to study time series obtained by the observation of real systems.
\end{abstract}

\section{Introduction} \label{sec:intro}

In this paper we are interested in studying the behaviour of the
Kolmogorov-Sinai (KS) entropy for randomly perturbed dynamical
systems. As it has been proved in \cite{kifer:book} for a typical
situation, the KS entropy of a dynamical system perturbed by a
sequence of random variables is infinite for all $\sigma >0$; here
$\sigma$ is a real parameter measuring the size of the
perturbation. This result shows that the KS entropy is not the right
quantity to look at in randomly perturbed dynamical systems. In
\cite{kifer:book} it has been suggested to study the KS entropy
relative to a partition.

We study the perturbation of a dynamical system $(X,\mu,f)$ by some
noise. By noise we mean a discrete stochastic processes of independent
and identically distributed (i.i.d.) random variables $\{ w_n \}_{n
\in \N}$ defined on a probability space $(W,{\cal P},P)$ with values
on the interval $[-\sigma,\sigma]$, for a positive real number
$\sigma$. The number $\sigma$ is the parameter that measures the
size of the perturbation on the original system.

We also assume that for all $\sigma$ there is an ergodic invariant
probability measure $\mu_\sigma$ on the perturbed system
$(X,f_\sigma)$ such that the probability measure $\mu$ is the weak
limit of the sequence $\mu_\sigma$ as $\sigma \to 0$ (for the formal
definition of $\mu_\sigma$ see Section \ref{sec:estim}).

In \cite{kifer:book} it has been proved that, under these hypothesis,
it holds
\begin{equation} \label{eq:kifer}
\limsup\limits_{\sigma \to 0^+} \ h_{\mu_\sigma} (f_\sigma, Z) \le
h_\mu (f,Z) \le h_\mu (f)
\end{equation}
where $h_\mu (f,Z)$ and $h_{\mu_\sigma} (f_\sigma, Z)$ are the KS
entropies relative to the partition $Z$ of the original and the
perturbed system respectively, and by $h_\mu (f)$ we denote the KS
entropy of the unperturbed dynamical system. The equality holds for
certain hyperbolic systems (\cite{kifer:book}). 

In equation (\ref{eq:kifer}) it is used the KS entropy relative to a
fixed partition, avoiding the problem of having infinite KS entropy,
and it is performed the limit for the size of the perturbation going
to zero. This result shows that for an unperturbed dynamical system,
the KS entropy relative to a partition is well approximated by the
same quantity for the perturbed system, when the perturbation is small
enough.

Unfortunately when we want to study time series obtained by
measurements of real systems, it is impossible to know how small is
the perturbation to the original system. Hence it is impossible to
perform the limit in equation (\ref{eq:kifer}). 

Real systems are observed using a symbolic representation of the
measurements. That is the phase space is divided into a finite number
of sets, and we study the time series given by the symbols of the sets
of the partition visited by the system as time grows.

The main idea to study real systems is to look at the system using
different scales of observation, which corresponds to use different
sizes for the partition of the phase space. When the partition is
coarse the effects of the noise are negligible, whereas, when the
partition is very fine, then the effects of the noise hide the
dynamics of the original system. Then the right approach to randomly
perturbed systems is to study the behaviour of the KS entropy while
the partition is refined and the noise is fixed.

In this paper we prove that, when the diameter of the partition is big
with respect to the size of the noise, the KS entropy of the perturbed
system is a good approximation for the unperturbed KS entropy. Our
main result is contained in equations (\ref{eq:e>s}) and
(\ref{eq:e<s}), which in particular imply equation
(\ref{eq:kifer}). Moreover we obtain precise estimates for the effects
of the noise, in order to distinguish between the perturbation and the
dynamics of the system, a very important problem in the analysis of
time series (\cite{cencini}). We remark that the main extension of the
result contained in equation (\ref{eq:kifer}) is its applicability to
time series obtained by observations of real systems.

To apply our method we have to face the problem of the measurement of
the KS entropy in real systems (that is the KS entropy of observed
time series). This is an important problem (for a review see for
example \cite{abarbanel}), and many papers have been dedicated to the
approximation of the KS entropy by different indicators. We recall the
approach of Grassberger and Procaccia \cite{grass-proc-k2} based on
the generalized Renyi entropy, and that of Cohen and Procaccia
\cite{cohen-proc-ks} based on the correlation integral for a time
series (\cite{grass-proc-corr-int}). In this paper we use a method
introduced in \cite{argenti} and \cite{licatone}, which is related to
the notion of {\it Algorithmic Information Content (AIC)} (see Section
\ref{sec:aic} for the definitions).

Moreover, when applying theoretical methods to real time series, it is
always questionable whether all the hypothesis needed for the methods
are verified by the single observations we study. This is this case,
for example, for ergodicity that we suppose to be verified.

In the following section, we recall some basic definitions and results
related to the KS entropy and its computation for dynamical systems
and purely random systems. In Sections \ref{sec:estim} and
\ref{sec:aic_noise} we apply the method of the AIC to randomly
perturbed systems, and finally in Section \ref{sec:num} we show some
numerical experiments performed on the logistic map at the chaotic
parameter $\lambda=4$.

\section{Basic definitions and results}

We now briefly recall the basic definitions we need. Let $(X,{\cal
B},\mu)$ be a probability space, where $X$ is a compact metric space
and $\cal B$ is the Borel $\sigma$~-algebra. In this paper we restrict
our attention on one dimensional spaces $X$, but we believe that the
techniques used can be generalized to higher dimensions. Let $f:X\to
X$ be a continuous map, invariant with respect to $\mu$ and
ergodic. Given a finite measurable partition $Z=\{I_i\}_{i=1,\dots,N}$
of the space $X$, the {\it entropy $H_\mu(Z)$} of the partition is
defined as
\begin{equation} \label{eq:ent-part}
H_\mu(Z)= - \sum_{i=1}^N \ \mu(I_i) \ \log (\mu(I_i))
\end{equation}
Let $f^{-1} Z$ be the partition given by the counter images $f^{-1}
I_i$. Then let 
\begin{equation} \label{eq:iter-part}
Z_n = Z\vee f^{-1} Z \vee \dots \vee f^{-n+1} Z
\end{equation}
be the partition given by the sets of the form 
$$I_{i_0} \cap f^{-1} I_{i_1} \cap \dots \cap f^{-n+1} I_{i_{n-1}}$$
varying $I_{i_j}$ among all the sets of $Z$. Then the {\it
Kolmogorov-Sinai entropy $h_\mu (f,Z)$ relative to the partition $Z$}
is defined as the limit
\begin{equation} \label{eq:ks-ent-part}
h_\mu (f,Z) = \lim\limits_{n\to \infty} \frac{1}{n} H_\mu (Z_n)
\end{equation}
The {\it Kolmogorov-Sinai entropy $h_\mu (f)$} of the dynamical system
$(X,\mu,f)$ is defined as
\begin{equation} \label{eq:ks-ent}
h_\mu (f) = \sup \{ h_\mu (f,Z) \ / \ Z \mbox{ finite partition} \} 
\end{equation}
Moreover for the explicit computation of the KS entropy, the
Kol\-mo\-go\-rov-Sinai theorem says that there are some special
partitions, called {\it generating}, for which the supremum among all
the partitions is realized. Hence it is enough to compute the KS
entropy relative to a generating partition to obtain the KS entropy of
the system. The existence of a generating partition for dynamical
systems is given by the following theorem:

\begin{teorema}[Krieger Generator Theorem \cite{krieger}] \label{teo:krieg}
Let $(X,\mu,f)$ be an ergodic dynamical system on a Lebesgue space
$X$, such that the probability measure $\mu$ is invariant and $h_\mu
(X,f) < \infty$. Then there is a finite generating partition for $f$.
\end{teorema}

Let $\epsilon$ be the diameter of a uniform partition $Z$ of the space
$X$, and simply denote by $h_\mu(\epsilon)$ the KS entropy $h_\mu
(f,Z)$ relative to the partition $Z$. This notation has the aim to
enhance the role of the diameter of the partition considered. The
quantity $h_\mu(\epsilon)$ is also called the {\it $\epsilon$-entropy}
of the dynamical system (see \cite{gw:noise} for a review). The
Krieger Generator Theorem implies that, if $h_\mu (f)$ is finite, there
exists $\epsilon_0$ such that $h_\mu (\epsilon)=h_\mu (f)$ for all
$\epsilon \le \epsilon_0$.

It is also possible to establish the behaviour of the
$\epsilon$-entropy for discrete stochastic processes of i.i.d. random
variables.

\begin{teorema}[Gaspard-Wang \cite{gw:noise}] \label{teo:gw}
Let $\{w_n \}_{n\in \N}$ be a discrete stochastic process of
i.i.d. random variables with values on the interval $[0,1]$. Then as
$\epsilon\to 0$ it holds $$h(\epsilon) \sim - \log \epsilon$$ where
the entropy is computed with respect to the invariant probability
measure of the system (the induced measure) absolutely continuous with
respect to the Lebesgue measure.
\end{teorema}

When dealing with randomly perturbed systems, it has been shown using
numerical experiments (\cite{gw:noise},\cite{deco}) that the expected
behaviour is
\begin{equation}\label{eq:stime}
h_{\mu_\sigma} (\epsilon) := h_{\mu_\sigma} (f_\sigma,Z) \sim \left\{
\begin{array}{lr}
h_\mu (\epsilon) & \hbox{for } \epsilon >> \sigma \\[2mm]
- \log \epsilon & \hbox{for } \epsilon << \sigma
\end{array}
\right.
\end{equation}
where $\sigma$ denotes the standard deviation of the noise (that is
the square root of the variance of the random variables $w_n$). The
same result is expected for the generalized Renyi entropy with $q=2$
(\cite{schreiber}).

We recall that for randomly perturbed systems the KS entropy
$h_{\mu_\sigma} (f_\sigma)$ is infinite, hence the Krieger Generator
Theorem is not applicable to the system $(X,f_\sigma)$. But we assume
that the unperturbed system $(X,f)$ has finite KS entropy $h_\mu (f)$
with generating partition of diameter $\epsilon_0$. Hence $h_\mu
(\epsilon) = h_\mu (f)$ for all $\epsilon \le \epsilon_0$. This
implies that if the size $\sigma$ of the noise is small with respect
to $\epsilon_0$, approximating $h_\mu (\epsilon)$ by $h_{\mu_\sigma}
(\epsilon)$ we find a good approximation to the KS entropy of the
unperturbed system.

In the following we prove (\ref{eq:stime}) using the {\it Algorithmic
Information Content (AIC)}, briefly described in the following
section. We remark that the behaviour of equation (\ref{eq:stime}) is
of great importance since it would give us a method to estimate the
size of the random perturbation on the deterministic dynamics of the
system we are observing.

\section{The algorithmic information content} \label{sec:aic}

Given a finite alphabet $\cal A$, let ${\cal A}^n$ be the set of all
words on the alphabet $\cal A$ of length $n$. The intuitive meaning of
quantity of information contained in a finite string $s \in {\cal
A}^n$ is the length of the smallest message from which you can
reconstruct $s$ on some machine. Thus, formally, the information $I$
is a function
$$ I:\mathcal{A}^{\ast } = \bigcup_{n\in \N} {\cal A}^n \rightarrow \N$$
on the set of finite strings on a finite alphabet
$\mathcal{A}$ which takes values in the set of natural numbers. 

One of the most important information function is the {\it Algorithmic
  Information Content ($AIC$)}. In order to define it, it is necessary
to define the notion of partial recursive function. We limit ourselves
to give an intuitive idea which is very close to the formal
definition. We can consider a partial recursive function as a computer
$C$ which takes a program $P$ (namely a binary string) as an input,
performs some computations, and gives a string $s=C(P)$, written on
the given alphabet $\mathcal{A}$, as an output. The $AIC$ of a string
$s$ is defined as the shortest binary program $P$ which gives $s$ as
its output, namely
\begin{equation} \label{eq:aic}
{AIC}(s,C)=\min \{|P|:C(P)=s\}
\end{equation}
We require that our computer is a universal computing machine. Roughly
speaking, a computing machine is called \emph{universal} if it can
simulate any other machine. In particular every real computer is a
universal computing machine, provided that we assume that it has
virtually infinite memory. For a precise definition see for example
\cite{livi} or \cite{cha}. We have the following theorem

\begin{teorema}[Kolmogorov \cite{kolmogorov}] \label{teo:aic-kolm}
If $C$ and $C^{\prime }$ are universal computing machines then
$$ \left| {AIC}(s,C)-{AIC}(s,C^{\prime })\right| \leq K\left(
C,C^{\prime }\right)$$ where $K\left( C,C^{\prime }\right) $ is a
constant which depends only on $C$ and $C^{\prime }$.
\end{teorema}

This theorem implies that the information content ${AIC}$ of $s$ with
respect to $C$ depends only on $s$ up to a fixed constant, then its
asymptotic behaviour does not depend on the choice of $C$. For this
reason from now on we will write ${AIC}(s)$ instead of ${AIC}(s,C)$.

The shortest program which gives a string as its output is a sort of
encoding of the string, and the information which is necessary to
reconstruct the string is contained in the program. Unfortunately the
coding procedure associated to the Algorithmic Information Content
cannot be performed by any algorithm. This is a very deep statement
and, in some sense, it is equivalent to the Turing halting problem or
to the G\"odel incompleteness theorem. Then the Algorithmic
Information Content is a function not computable by any algorithm.

Using the notion of Algorithmic Information Content it is possible to
define a notion of {\it complexity} for infinite strings. Let $\omega$
be an infinite string on the alphabet $\cal A$, that is $\omega \in
\Omega := {\cal A}^\N$. We denote by $\omega^n$ the first $n$ symbols
of the string $\omega$. Then $\omega^n \in {\cal A}^n$. The complexity
measures the mean quantity of information in each digit of the string
$\omega$. Formally

\begin{definizione}[Brudno \cite{brudno}] \label{def:compl}
The {\it complexity $K(\omega)$} of an infinite string $\omega \in
\Omega$ is given by $$K(\omega) = \limsup\limits_{n\to \infty} \
\frac{AIC(\omega^n)}{n}$$
\end{definizione}

Using the method of symbolic dynamics it is possible to consider the
information and the complexity of the orbits of a dynamical
system. Let $(X,\mu,f)$ be an ergodic dynamical system and let
$Z=\{I_1,\dots,I_N \}$ be a finite measurable partition of the space
$X$. To the partition $Z$ it is associated the finite alphabet ${\cal
A}=\{1,\dots,N \}$. Define a map $\varphi_Z: X \to \Omega$, where
$\Omega = {\cal A}^{\N}$, by
\begin{equation}\label{eq:orb-simb}
(\varphi_Z (x))_j = k \quad \Longleftrightarrow \quad f^j(x) \in I_k
\end{equation}
The sequence $\varphi_Z(x)\in \Omega$ associated to a point $x\in X$
is called the {\it symbolic orbit} of $x$ relative to the partition
$Z$. Then we have the following definition.

\begin{definizione}[Brudno \cite{brudno}] \label{def:compl-ds}
The {\it complexity $K(x,Z)$ relative to the partition $Z$} of the
orbit with initial condition $x\in X$ is given by
$$K(x,Z) := K(\varphi_Z (x)) = \limsup\limits_{n\to \infty} \
\frac{AIC((\varphi_Z (x))^n)}{n}$$
\end{definizione}

\begin{guess} \label{rem:1}
  In \cite{brudno} and in \cite{gal-cs}, using open covers and
  computable structures, notions of complexity of orbits of a
  dynamical system are defined that do not depend on the partition.
  But in this paper we are only interested in the complexity dependent
  on a partition for the same reason for which we study KS entropy
  relative to a partition.
\end{guess}

This notion of complexity has been related to the notion of KS entropy
by the following theorem.

\begin{teorema}[Brudno \cite{brudno}] \label{teo:brudno}
Let $(X,\mu,f)$ be an ergodic dynamical system and $Z$ be a finite
measurable partition, then $K(x,Z)=h_\mu (f,Z)$
for $\mu$-almost any $x\in X$.
\end{teorema}

\begin{teorema}[White \cite{white1},\cite{white2}] \label{teo:white}
In the same hypothesis as before, for $\mu$-almost any $x\in X$ it
holds
$$\liminf\limits_{n\to \infty} \ \frac{AIC((\varphi_Z (x))^n)}{n} =
h_\mu (f,Z)$$
\end{teorema}

If the AIC were a computable function, using the previous theorems we
could compute the KS entropy relative to a partition. This approach
can still be useful using {\it optimal compression algorithms}, that
is algorithms which encode symbolic strings, giving an approximation
of the information contained in a string, hence an approximation of
its AIC. For this approach see \cite{licatone}, where also a new
compression algorithm is presented and applied to some well known
chaotic systems.

\begin{guess} \label{rem:2}
The notion of complexity for dynamical systems has been studied for
some well known weakly chaotic systems (\cite{bon-men},\cite{bona})
for which it holds $h_\mu(f,Z)=K(x,Z)=0$ for all the partitions, with
the aim of giving a classification of these systems according to the
asymptotic behaviour of the AIC of the orbits.
\end{guess}

\begin{guess} \label{rem:3}
The notion of AIC has been linked to other indices for dynamical
systems. For example it has been related to the sensitivity to initial
conditions (\cite{gal}) and to the Poincar\'e recurrence times
(\cite{bgi}).
\end{guess}

\section{Estimates of the relative KS entropy} \label{sec:estim}

We now introduce formal definitions for random perturbations of a
dynamical system (following \cite{kifer:book}). We also give some
estimates for the KS entropy relative to a partition for the perturbed
system using the tool of symbolic dynamics.

Let $\{ w_n \}$ be a stochastic process of i.i.d. random variables,
with each $w_n$ defined on the probability space $(W,{\cal P},P)$ with
values on the interval $[-\sigma,\sigma]$, and $q^\sigma$ as the
induced distribution.

\begin{definizione} \label{def:rand-pert}
A {\it random perturbation} of the dynamical system $(X,\mu,f)$ is a
family of Markov chains $f^n_\sigma$ on the probability space
$(W,{\cal P},P)$ with values on $X$ and with transition probabilities
given by
$$P^\sigma (x,B) := P \{ f^{n+1}_\sigma \in B \ |\ f^n_\sigma = x \} =
q^\sigma (B-f(x))$$ for all Borel sets $B\subset X$.
\end{definizione}

The meaning of this definition is that the point $x$ is moved to the
point $f(x)$ under the unperturbed dynamics and then it disperses
randomly with distribution $q^\sigma$. 

We will consider two different possible actions of the random
perturbation (see \cite{argyris}). We talk of {\it output noise} when
the dynamics is driven only by the unperturbed map $f$, and the
dispersion is caused by a non exact observation of the point
$f(x)$. That is given a point $x\in X$, its orbit is given by
$x_n=f(x_{n-1})$ and our data are measured observing the orbit $y_n =
f_\sigma (x_n) = x_n + w_n$. Instead we talk of {\it dynamical noise}
when the dispersion is intrinsic in the dynamics and it is not caused
by the observation, in this case the dynamics is driven by
$f_\sigma$. So, for a point $x\in X$ we have $x_n = f^n_\sigma (x) =
f(x_{n-1}) + w_n$.

\begin{definizione} \label{def:rand-meas}
A probability measure $\mu_\sigma$ on $X$ is an {\it invariant
measure} of the randomly perturbed dynamical system $(X,f_\sigma)$ if
$$\mu_\sigma (B) = \int_X P^\sigma (x,B) \ d \mu_\sigma (x)$$
for all Borel sets $B\subset X$.
\end{definizione}

We now start from the following assumptions:
\begin{description}
\item[(i)] there is an ergodic dynamical system $(X,\mu,f)$ with
finite KS entropy, and denote by $h_\mu (\epsilon)$ its KS entropy
relative to a finite measurable partition $Z$ of diameter $\epsilon$;
\item[(ii)] we can analyze data produced by a random perturbation
$(X,f_\sigma)$ of the system (output or dynamical noise) with fixed
$\sigma$;
\item[(iii)] there is an invariant and ergodic probability measure
$\mu_\sigma$ on $(X,f_\sigma)$ and the KS entropy relative to the
partition $Z$ for the randomly perturbed system is denoted by
$h_{\mu_\sigma} (\epsilon)$.
\end{description}

To analyze our system by the Algorithmic Information Content we have
to use the method of symbolic dynamics. Given the finite measurable
partition $Z=\{I_1,\dots,I_N \}$ of diameter $\epsilon = \frac{1}{N}$,
we can associate a symbolic orbit in $\Omega = \{1,\dots,N\}^\N$ to
any orbit in the space $X$ using the map $\varphi_Z$ defined in
equation (\ref{eq:orb-simb})

In the case of the unperturbed system $(X,\mu,f)$ we denote by
$\Omega' \subseteq \Omega$ the image of the map $\varphi_Z$. The set
$\Omega'$ is closed with respect to the action of the usual shift map
$\tau$ defined on $\Omega$. It is also possible to define the
probability measure $\nu$ on $\Omega'$ induced by $\varphi_Z$. So we
can work on the ergodic dynamical system $(\Omega',\nu,\tau)$ with KS
entropy $h_\nu = h_\mu (\epsilon)$.

Analogously we define the map $\varphi_{Z,\sigma}$ from the product
space $X\times W$ to $\Omega$ in the following way: for each $x\in X$
the map $\varphi_{Z,\sigma}(x,\cdot)$ is the symbolic map for the
Markov chain $f^n_\sigma$. So for every point $x\in X$ there is a set
of possible symbolic orbits in $\Omega'_\sigma :=\varphi_{Z,\sigma}
(X\times W)$ that correspond to different realizations of the
stochastic process $\{ w_n \}$. Let $\nu_\sigma$ be the probability
measure induced on $\Omega'_\sigma$ by the measure $\mu_\sigma$. Then
we study the ergodic dynamical system
$(\Omega'_\sigma,\nu_\sigma,\tau)$ with KS entropy $h_{\nu_\sigma} =
h_{\mu_\sigma} (\epsilon)$.

Since now we will denote by $\psi$ and by $\omega$ the symbolic orbits
in $\Omega'$ and $\Omega'_\sigma$ respectively. 

From the previous definitions, the following propositions easily follow 

\begin{prop} \label{prop:spazi-simb}
$\Omega' \subseteq \Omega'_\sigma$
\end{prop}

\begin{prop} \label{prop:ent-simb}
$h_\mu (\epsilon) \le - \log \epsilon \ , \quad h_{\mu_\sigma}
(\epsilon)\le - \log \epsilon \ .$
\end{prop}

\begin{prop} \label{prop:spazi-ug}
If $\Omega' = \Omega'_\sigma$ and $\nu$ is absolutely continuous with
respect to $\nu_\sigma$ ($\nu << \nu_\sigma$) then $h_\mu (\epsilon) =
h_{\mu_\sigma}(\epsilon)$.
\end{prop}

\noindent {\bf Proof.} For $\nu_\sigma$-almost any $\omega \in
\Omega'_\sigma$ it holds
\begin{equation} \label{eq:si-rum}
\lim\limits_{n\to \infty} \ \frac{AIC(\omega^n)}{n} =
h_{\mu_\sigma}(\epsilon)
\end{equation}
thanks to Theorems \ref{teo:brudno} and \ref{teo:white}. Hence, since for
$\nu$-almost all $\psi \in \Omega'$ it holds
\begin{equation} \label{eq:no-rum}
\lim\limits_{n\to \infty} \ \frac{AIC(\psi^n)}{n} =
h_{\mu}(\epsilon)
\end{equation}
and $\nu << \nu_\sigma$ then $h_\mu (\epsilon) =
h_{\mu_\sigma}(\epsilon)$. \qed

\vskip 0.5cm
We conclude this section with a lower bound for the $\epsilon$-entropy
of the random perturbed system. 

\begin{prop}[Kifer \cite{kifer:book} Theorem 2.4] \label{teo:kifer}
Suppose that all transition probabilities $P^\sigma(x,\cdot)$ have
bounded densities $p^\sigma(x,y)\le K$, then $h_{\mu_\sigma}(\epsilon)
\ge -\log \epsilon - \log K$.
\end{prop}

Hence the KS entropy of a randomly perturbed system is always
infinite, since it can be computed as the limit of $h_{\mu_\sigma}
(\epsilon)$ for $\epsilon$ going to 0.

\section{The AIC for noisy systems} \label{sec:aic_noise}

\subsection{The case of output noise}
\label{sec:output} 

We now study the behaviour of the AIC relative to a partition in
randomly perturbed dynamical systems. We start with the case of output
noise (see previous section) which is easier to be studied and hence
it is useful to introduce the arguments we will use for the more
important case of dynamical noise.

Let $\Omega'$ and $\Omega'_\sigma$ be the symbolic spaces of the
unperturbed and the perturbed systems respectively. We will study the
$\epsilon$-entropy of the perturbed system using the AIC of the
symbolic orbits, $\psi$ and $\omega$, of the original system and of
its perturbed version.

At any step of our dynamical system the perturbation induced by the
noise changes the position of the point, and so it could change the
set of the partition in which the point is. So the $n$-th symbol of
the string $\omega$ could be different from $\psi_n$, the
corresponding symbol of the string $\psi$. We recall that since we are
in the case of output noise, at the step $n+1$ there is no memory of
the noise at the previous step. Moreover since the random variables
$w_n$ are independent, the probability $p$ of $\omega_n$ being
different from $\psi_n$ does not depend on the step $n$, neither on
the set of the partition we consider. 

Let $H(p)$ be the KS entropy of a {\it $p$-Bernoulli trial}, namely
the stochastic process of independent Bernoulli variables
$\vartheta_n$ with parameter $p$. It is well known that
$$H(p)=-p\ \log (p) - (1-p)\ \log (1-p) \ .$$ Hence applying
Theorems \ref{teo:brudno} and \ref{teo:white} it holds
\begin{equation} \label{eq:p_ber}
\lim\limits_{n\to \infty} \frac{H_n(p)}{n} = H(p)
\end{equation}
where $H_n(p)$ denotes the average of the AIC over all the
sequences produced by the $p$-Bernoulli trial.

In this framework, we have

\begin{prop} \label{prop:disug-alta}
If the random variables $w_n$ have values in the interval
$[-\sigma,\sigma]$, then
$$\E_{\mu_\sigma} [ AIC(\omega^n)] \leq \E_\mu [AIC(\psi^n)] + np \log
\left(2 \left\lceil \frac{\sigma}{\epsilon} \right\rceil \right) +
H_n(p)$$ where $\lceil \cdot \rceil$ denotes the upper integer
part of a real number and $\E[\cdot]$ denotes the mean value.
\end{prop}

\noindent {\bf Proof.} To prove the proposition we show an algorithm
to describe the string $\omega^n$ and compute the information it
needs. Let assume that we know the string $\psi^n$. We have then
just to specify the symbols of $\omega^n$ that are different from those
of $\psi^n$. We call $i_n$ this information. We show that for all
strings $\omega^n$ and $\psi^n$ it holds
\begin{equation}\label{eq:r_n}
i_n \leq r_n = \left( \# \{ i\ /\ \psi_i \not= \omega_i \} \right) \
\log \left(2 \left\lceil \frac{\sigma}{\epsilon} \right\rceil \right)
+ H_n(p)
\end{equation} 
The algorithm to describe $\omega^n$ works as follows. First it
describes a binary string $s=(s_0,\dots,s_{n-1})$, where $s_i=0$
implies that $\omega_i=\psi_i$ and $s_i=1$ otherwise. To this aim we
need on average $H_n(p)$ bits of information, since the string $s$ is
obtained as a string of a $p$-Bernoulli trial.

Moreover the algorithm needs also to specify how far on the left or on
the right the noise has moved the point. To this aim we need $\log
\left(2 \left\lceil \frac{\sigma}{\epsilon} \right\rceil \right)$ bits
of information for each symbol of $\omega^n$ different from the
corresponding symbol of $\psi^n$. The factor $\left\lceil
\frac{\sigma}{\epsilon} \right\rceil$ counts exactly how many sets of
the partition $Z$ can be covered by the effect of the noise, and the
factor 2 is needed to specify whether the point is moved to the right
or to the left.

Hence to specify the string $\omega^n$ we do not need more than
$AIC(\psi^n)+r_n$ bits of information. If we evaluate the mean of
$AIC(\omega^n)$ over the measure $\mu_\sigma$ then the thesis follows
from the definition of the probability $p$. \qed

\begin{teorema} \label{teo:stima-noise}
If $h_\mu (\epsilon)$ is the KS entropy relative to a partition $Z$ of
diameter $\epsilon$ of an unperturbed dynamical system $(X,\mu,f)$,
for the $\epsilon$-entropy $h_{\mu_\sigma} (\epsilon)$ of the system
perturbed by an output noise $\{ w_n \}$ with values in the interval
$[-\sigma,\sigma]$ it holds
$$
h_{\mu_\sigma} (\epsilon) \leq h_\mu (\epsilon) + p \log \left(2
\left\lceil \frac{\sigma}{\epsilon} \right\rceil \right) + H(p) 
$$
\end{teorema}

\noindent {\bf Proof.} The thesis follows from Proposition
\ref{prop:disug-alta}, using equations (\ref{eq:no-rum}) and
(\ref{eq:si-rum}) for $h_\mu (\epsilon)$ and $h_{\mu_\sigma}
(\epsilon)$, and using equation (\ref{eq:p_ber}). \qed

\subsection{The case of dynamical noise} 
\label{sec:dynamical}

We now study the case of dynamical noise. We assume to have an
unperturbed dynamical system $(X,\mu,f)$ and a dynamical perturbation
due to a stochastic process $\{ w_n \}$ of i.i.d. random variables
with values in the interval $[-\sigma,\sigma]$, such that the points
$x\in X$ follow orbits given by $x_{n+1}= f(x_n) + w_{n+1}$.

Let consider a partition $Z$ of the space $X$ in a finite number $N$
of measurable sets $\{I_i \}$, such that the diameter of each set of
the partition is $\epsilon = 1/N$. Let $\Omega=\{1,\dots,N\}^\N$ be
the space of the symbolic orbits associated to the dynamical system,
and denote by $\omega$ a single string in $\Omega'_\sigma$, given by
the measurements made on the system. That is $\omega$ is the symbolic
orbit of a point $x\in X$ after the effect of the noise. To estimate
the Algorithmic Information Content of the string $\omega$ from above,
we show how to construct it from the knowledge of the noise and of the
unperturbed dynamics.

Let $x_0$ denote the initial point of an orbit of the system. By the
classical symbolic representation of orbits, there is a symbolic
string $\psi^0=(\psi^0_0,\psi^0_1,\dots) \in \Omega$, relative to the
partition $Z$, associated to the point $x_0$ for the unperturbed
dynamical system, that is $f^i(x_0)\in I_{\psi^0_i}$ for all $i\in\N$,
and there is also the symbolic string
$\omega=(\omega_0,\omega_1,\dots)$ that is given by the perturbed
system. We have $\omega_0=\psi^0_0$. In the case of dynamical noise,
we have to consider the effect due to the fact that the action of the
noise is not forgotten at each step. So we have to extend the
knowledge of the unperturbed orbit to more than just one iteration. To
find a good number of iterations to be stored, we use the
approximation of the Kolmogorov-Sinai entropy
$h_\mu(f,Z)=h_\mu(\epsilon)$ relative to the partition $Z$ by the
decreasing sequence $H(T^{-n}Z | Z_n)$, where $Z_n= Z\vee T^{-1}Z \vee
\dots \vee T^{-n+1}Z$ and
$$H(P|Q)=-\sum_{i,j} \ \mu(P_i\cap Q_j) \ \log \frac{\mu(P_i\cap
Q_j)}{\mu(Q_j)}$$
for any two finite partitions $P=(P_i)$ and $Q=(Q_j)$. Hence for any
$\delta >0$ there exists $n_0 \in \N$ such that 
\begin{equation} \label{eq:approx-ent}
h_\mu(\epsilon)- \delta \le H(T^{-n_0} Z | Z_{n_0}) \le
h_\mu(\epsilon) + \delta
\end{equation}

Let $\delta>0$ be fixed, and find the corresponding integer $n_0$
according to equation (\ref{eq:approx-ent}). Assume that for the
initial point $x_0$ of the orbit, we know the first $n_0$ digits of
the sequence $\psi^0$. That is we assume to know the cylinder
$Z_{n_0}(x_0)$, where with this notation we mean the set of the
partition $Z_{n_0}$ containing $x_0$.

Let $\tilde x_1=f(x_0)$ and $x_1 =\tilde x_1 + w_1$. As before we
denote by $\psi^1$ the unperturbed symbolic orbit associated to the
point $x_1$. We assume to know the action of the noise $w_1$ with
respect to the partition $Z_{n_0}$, so we compare the strength of the
noise $\sigma$ with the diameter\footnote{In this case we cannot ask
the partition $Z_{n_0}$ to be uniform, so we consider the diameter to
be the smaller diameter of the intervals that make the partition
$Z_{n_0}$.} $\epsilon_{n_0}$ of the partition $Z_{n_0}$. Moreover, to
know $Z_{n_0}(\tilde x_1)$ we only need the symbol $\psi^0_{n_0}$, and
to know $Z_{n_0}(x_1)$, once we have $Z_{n_0}(\tilde x_1)$, we only
need the action of $w_1$ with respect to the partition $Z_{n_0}$. This
is enough to obtain $\omega_1 = \psi^1_0$. So to obtain the first two
symbols of the string $\omega$ it is enough the following information:
$I(\psi^0_0,\dots, \psi^0_{n_0-1})$, $I(w_1|Z_{n_0})$ and
$I(\psi^0_{n_0}|\psi^0_0,\dots, \psi^0_{n_0-1})$. We have used
$I(\cdot | \cdot)$ to denote the conditional information.

At this point it is possible to iterate the argument with the string
$\psi^1$ instead of the string $\psi^0$, hence to find the symbol
$\omega_2 = \psi^2_0$ using the following information:
$I(w_2|Z_{n_0})$, $I(\psi^1_0,\dots, \psi^1_{n_0-1})$ (which is given
by $Z_{n_0}(x_1)$ and is indeed known) and
$I(\psi^1_{n_0}|\psi^1_0,\dots, \psi^1_{n_0-1})$.

Iterating this argument we obtain the following proposition:

\begin{prop} \label{prop:dyn-disug-alta}
If the random variables $w_n$ have values in the interval
$[-\sigma,\sigma]$, then for any $\delta >0$ there exists $n_0 \in \N$
and a constant $C>0$ such that
$$\E_{\mu_\sigma} [ AIC(\omega^n)] \leq (n-1) \left[ h_\mu(\epsilon) +
\delta + p \log \left(2 \left\lceil \frac{\sigma}{\epsilon_{n_0}}
\right\rceil \right) \right] + H_{n-1}(p)+C$$ where $\lceil \cdot
\rceil$ denotes the upper integer part of a real number, the
probability $p$ is defined as in Proposition \ref{prop:disug-alta} and
$\omega^n = (\omega_0, \dots, \omega_{n-1})$.
\end{prop}

\noindent {\bf Proof.} To obtain the entropy on the right hand side we
only have to use equation (\ref{eq:approx-ent}), and to remark that
the conditional entropy $H(T^{-n_0}Z | Z_{n_0})$ is obtained as the
integral over the space of all admissible strings of the information
function $I(\psi_{n_0}|\psi_0,\dots, \psi_{n_0-1})$. The constant $C$
is the information $I(\psi^0_0,\dots, \psi^0_{n_0-1})$ needed to start
the iteration.

For the part related to the noise we repeat the same argument of
Proposition \ref{prop:disug-alta}. Then the action of the noise can be
bounded by
$$\E_P [r_n]=(n-1)p \log \left(2 \left\lceil
\frac{\sigma}{\epsilon_{n_0}} \right\rceil \right) + H_{n-1}(p)$$
bits of information. \qed

\vskip 0.5cm Let again $h_{\mu_\sigma} (\epsilon)$ be the
$\epsilon$-entropy of the perturbed dynamical system, hence it follows

\begin{teorema} \label{teo:dyn-stima-noise}
If $h_\mu (\epsilon)$ is the $\epsilon$-entropy of an unperturbed
dynamical system $(X,\mu,f)$, for the $\epsilon$-entropy
$h_{\mu_\sigma} (\epsilon)$ of the system perturbed by a dynamical
noise $\{ w_n \}$ it holds: for any $\delta >0$ there exists $n_0 \in
\N$ such that
$$
h_{\mu_\sigma} (\epsilon) \leq h_\mu (\epsilon) + \delta + p \log
\left(2 \left\lceil \frac{\sigma}{\epsilon_{n_0}} \right\rceil \right)
+ H(p) 
$$
\end{teorema}

\subsection{Conclusions} \label{sec:concl}

At this point we can make some conclusions about the estimates of the
KS entropy of a randomly perturbed system.

We recall the conditions we imposed in Section \ref{sec:estim} on the
data we study. So, given a fixed $\sigma$, we vary the diameter
$\epsilon$ of the partition considered.

Let $\epsilon > \sigma$. In this case we can assume that for some
$\epsilon$ we have $\Omega'=\Omega$. This implies (Propositions
\ref{prop:spazi-simb} and \ref{prop:spazi-ug}) that $h_{\mu_\sigma}
(\epsilon) = h_\mu (\epsilon)$. When $\Omega' \subset \Omega$ we have
that, fixed $\delta>0$ such that $h_\mu (\epsilon) + \delta < -\log
\epsilon$, if $\epsilon_{n_0} > \sigma$
\begin{equation} \label{eq:e>s}
-\log \epsilon -const \ \le h_{\mu_\sigma}(\epsilon) \le h_\mu
 (\epsilon) + \delta + p \log 2 + H(p) 
\end{equation}
using Proposition \ref{teo:kifer} and Theorem
\ref{teo:dyn-stima-noise}. We recall that $H(p)$ takes its maximum
value $\log 2$ for $p=0.5$, and it converges to $0$ as $p$ goes to $0$
or $1$. Then if $\epsilon_{n_0} >> \sigma$ we have that $p$ converges
to $0$, hence $(p\log 2 + H(p))$ converges to $0$, and equation
(\ref{eq:e>s}) becomes $h_{\mu_\sigma} (\epsilon) \le h_\mu (\epsilon)
+ \delta$. In particular this implies equation (\ref{eq:kifer}). If
instead $\epsilon_{n_0}\le \sigma$ then we have the same conclusion as
in the following case.

If $\epsilon < \sigma$ then it holds 
\begin{equation} \label{eq:e<s}
\begin{array}{c}
-\log \epsilon -const \ \le h_{\mu_\sigma}(\epsilon) \le \\[0.4cm] \le
\min \left( -\log \epsilon  , h_\mu (\epsilon) + \delta + p \log
\left(2 \left\lceil \frac{\sigma}{\epsilon_{n_0}} \right\rceil \right)
+ H(p) \right) 
\end{array}
\end{equation}
using Propositions \ref{prop:ent-simb}, \ref{teo:kifer} and Theorem
\ref{teo:dyn-stima-noise}. Again if $\epsilon_{n_0} << \sigma$ then $p$
converges to $1$, hence $H(p)$ converges to $0$, and the estimate on
the right hand side of equation (\ref{eq:e<s}) becomes $h_\mu
(\epsilon) + \delta - \log \epsilon_{n_0} + const.$

From equations (\ref{eq:e>s}) and (\ref{eq:e<s}) it emerges that by
computing the KS entropy relative to partitions, we obtain an
indication of the smallness of the random perturbation present in the
data we are studying. Indeed if in the curve $h_{\mu_\sigma}
(\epsilon)$ we find an almost flat region ($\epsilon \ge \epsilon_1$),
and as $\epsilon$ decreases there is a change in the behaviour of the
curve, that almost behaves as $-\log(\epsilon)$ ($\epsilon\le
\epsilon_2$), then the interval $(\epsilon_2, \epsilon_1)$ should be a
good approximation for the size $\sigma$ of the perturbation. In the
next section we apply this method to a perturbed dynamical system,
using a compression algorithm to estimate the Algorithmic Information
Content.

Moreover this method can also be thought of as a method to distinguish
between purely stochastic systems and randomly perturbed dynamical
systems. Indeed for randomly perturbed systems equations
(\ref{eq:e>s}) and (\ref{eq:e<s}) suggest the presence of a change of
behaviour for the $\epsilon$-entropy at the level of the
perturbation. This change should not appear in purely stochastic
systems. However this aspect has to be examined further.

\section{Numerical experiments} \label{sec:num}

In the previous section we have proved that the Algorithmic
Information Content can be used for randomly perturbed dynamical
systems to obtain information on the size of the random
perturbation. We remark that unfortunately the AIC, as defined in
equation (\ref{eq:aic}), is not a computable function. So to compute
the information content of a string one needs to use a {\it
compression algorithm}, that is an algorithm which approximates the
value of the information function AIC.

Formally, a {\it compression algorithm} is defined as a recursive
reversible coding procedure $Z:{\cal A}^* \to \{0,1\}^*$ (for example,
the data compression algorithms that are in any personal
computer). The information content $I_Z$ of a finite string $s$
computed by the algorithm $Z$ is given by the binary length of the
compressed string, that is $I_Z(s)=|Z(s)|$.

Not all compression algorithms approximate the AIC. In \cite{licatone}
this point is discussed in details. For the purposes of this paper, it
is sufficient that 
$$\limsup\limits_{n\to \infty} \frac{I_Z(\omega^n)}{n} = K(\omega)$$
where $K(\omega)$ is the complexity of an infinite string $\omega$
(see Definition \ref{def:compl}). The algorithms with this property
are called {\it optimal}. 

In \cite{licatone} it is also introduced a compression algorithm
called {\it CASToRe}, which has been created with the specific aim of
studying weakly chaotic dynamical systems, namely chaotic systems
with null KS entropy.

Using CASToRe we have analyzed an example of randomly perturbed
dynamical system. We have chosen the logistic map
$$f:x\mapsto \lambda x(1-x)$$ on the interval $[0,1]$ with
$\lambda=4$. To this map we have added a dynamical noise given by
independent random variables $\{w_n\}$ uniformly distributed on the
interval $[-\sigma,\sigma]$.

Keeping $\sigma$ fixed we have varied the diameter of the partition
from 0.5 to 0.004=1/250, and we have computed the complexity of the
$10^6$-long symbolic orbits of the perturbed system relative to the
different partitions.

\begin{figure}[ht]
\psfig{figure=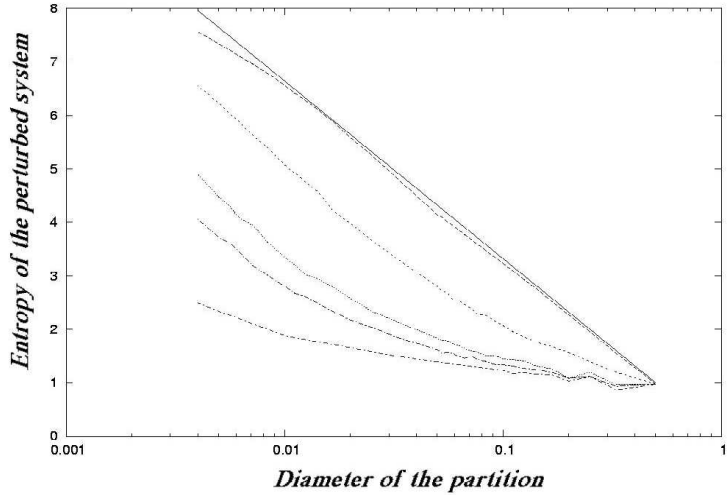,width=10cm,angle=270}
\caption{\it The complexity for different values of $\sigma$ with
respect to the diameter of the partition.} 
\label{fig:1}
\end{figure}

In figure \ref{fig:1} we show the results. It is a log-linear plot of
the complexity of the symbolic orbits versus the diameter of the
partition, plotted on the x axis. The solid line is the upper bound
$-\log \epsilon$, and the other curves are the empirical curves which
correspond to the following values of $\sigma$: 0.5, 0.1, 0.02, 0.01,
0.001 from the upper to the lower, respectively. For each of the
curves $h_{\mu_\sigma} (\epsilon)$ it is possible to identify the
intervals $(\epsilon_2,\epsilon_1)$, that give good approximation of
$\sigma$, in agreement with the theoretical result in equations
(\ref{eq:e>s}) and (\ref{eq:e<s}) and the subsequent comments.

When the diameter of the partition is high with respect to the size of
the noise, then we have a good approximation of the unperturbed KS
entropy of the system, since the partition with diameter $0.5$ is
generating. In figure \ref{fig:1} it is possible to see that the
complexity of the symbolic orbits for partitions with diameter close
to $0.5$ is close to $1$ (the KS entropy\footnote{The KS entropy is
$\log 2$, but we use logarithms in base $2$ since we measure binary
information contents.} of the logistic map) for the empirical curves
with $\sigma< 0.1$.


\begin{thebibliography}{10}

\bibitem{abarbanel} H.D.I.Abarbanel, R.Brown, J.J.Sidorowich,
L.S.Tsimring, {\it The analysis of observed chaotic data in physical
systems}, Rev. Mod. Phys. {\bf 65} (1993), 1331--1392

\bibitem{argenti} F.Argenti, V.Benci, P.Cerrai, A.Cordelli,
S.Galatolo, G.Menconi, {\it Information and dynamical systems: a
concrete measurement on sporadic dynamics}, Chaos Solitons Fractals
{\bf 13} (2002), 461--469

\bibitem{argyris} J.Argyris, I.Andreadis, {\it On the influence of
noise on the largest Lyapunov exponent of attractors of stochastic
dynamic systems}, Chaos Solitons Fractals {\bf 9} (1998),
959--963

\bibitem{licatone} V.Benci, C.Bonanno, S.Galatolo, G.Menconi,
M.Virgilio, {\it Dynamical systems and computable information},
Disc. Cont. Dyn. Syst. - B, to appear

\bibitem{bona} C.Bonanno, {\it The Manneville map: topological, metric
and algorithmic entropy}, http://arXiv.org/abs/math.DS/0107195 (2001)

\bibitem{bgi} C.Bonanno, S.Galatolo, S.Isola, {\it Poincar\'e
recurrence times and algorithmic complexity}, submitted (2003)

\bibitem{bon-men} C.Bonanno, G.Menconi, {\it Computational information for
the logistic map at the chaos threshold}, Disc. Cont. Dyn. Syst.- B
{\bf 2} (2002), 415--431

\bibitem{brudno} A.A.Brudno, {\it Entropy and the complexity of the
trajectories of a dynamical system}, Trans. Moscow Math. Soc. {\bf 2}
(1983), 127--151

\bibitem{cencini} M.Cencini, M.Falcioni, E.Olbrich, H.Kantz,
A.Vulpiani, {\it Chaos or noise: difficulties of a distinction},
Phys. Rev. E {\bf 62} (2000), 427--437

\bibitem{cha} G.J.Chaitin, ``Information, randomness and
incompleteness, papers on algorithmic information theory'', World
Scientific, Singapore, 1987

\bibitem{cohen-proc-ks} A.Cohen, I.Procaccia, {\it Computing the
Kolmogorov entropy from time signals of dissipative and conservative
dynamical systems}, Phys. Rev. A {\bf 31} (1985), 1872--1882

\bibitem{gal-cs} S.Galatolo, {\it Orbit complexity by computable
structures}, Nonlinearity {\bf 13} (2000), 1531--1546

\bibitem{gal} S. Galatolo, {\it Complexity, initial data sensitivity,
    dimension and weak chaos in dynamical systems}, Nonlinearity {\bf
    16} (2003), 1219--1238

\bibitem{gw:noise} P.Gaspard, X.J.Wang, {\it Noise, chaos, and
{$(\epsilon,\tau)$}-entropy per unit time}, Phys. Rep. {\bf 235}
(1993), 291--343 

\bibitem{grass-proc-k2} P.Grassberger, I.Procaccia, {\it Estimation of
the Kolmogorov entropy from a chaotic signal}, Phys. Rev. A {\bf 28}
(1983), 2591--2593

\bibitem{grass-proc-corr-int} P. Grassberger, I. Procaccia, {\it
Characterization of strange attractors}, Phys. Rev. Lett. {\bf 50}
(1983), 346--349

\bibitem{kifer:book} Y.Kifer, ``Random perturbations of dynamical
systems'', Birkh\"auser, 1988

\bibitem{kolmogorov} A.N.Kolmogorov, {\it Combinatorial foundations of
information theory and the calculus of probabilities},
Russ. Math. Surv. {\bf 38} (1983), 29--40

\bibitem{krieger} W.Krieger, {\it On entropy and generators of
measure-preserving transformations}, Trans. Amer. Math. Soc. {\bf
149} (1970), 453--464

\bibitem{livi} M.Li, P.Vitanyi, ``An introduction to Kolmogorov
complexity and its applications'', Springer, 1993

\bibitem{deco} C.Schittenkopf, G.Deco, {\it Identification of
deterministic chaos by an information theoretic measure of the
sensitive dependence on the initial conditions}, Physica D {\bf 110}
(1997), 173--181

\bibitem{schreiber} T.Schreiber, H.Kantz, {\it Noise in chaotic data:
diagnosis and treatment}, Chaos {\bf 5} (1995), 133--142

\bibitem{white1} H.White, ``On the algorithmic complexity of
trajectories of points in dynamical systems'', Ph.D. dissertation
Univ. of North Carolina at Chapel Hill, 1991

\bibitem{white2} H.White, {\it Algorithmic complexity of points in a
dynamical system}, Ergod. Th. Dynam. Syst. {\bf 13} (1993), 807--830

\end{thebibliography}
\end{document}